\documentclass[12pt]{article}
\usepackage{amsmath,amstext,amscd}
\usepackage{amssymb}
\usepackage{epsfig}
 \usepackage[amsmath,thmmarks]{ntheorem}
\usepackage{amsfonts,latexsym}
 \input epsf.tex

\theoremseparator{.}

{\theoremnumbering{Alph} \newtheorem{MainThm}{Theorem}}
\newtheorem{MainCor}{Corollary}[MainThm]

\newtheorem{thm}{Theorem}[section]
\newtheorem{lemma}[thm]{Lemma}
\newtheorem{cor}[thm]{Corollary}
\newtheorem{prop}[thm]{Proposition}

\theorembodyfont{\upshape}
\newtheorem{defn}[thm]{Definition}

\newtheorem{Remark}[thm]{Remark}

\newtheorem{Open questions}[thm]{Open questions}

\theoremheaderfont{\em}
\theoremstyle{nonumberplain}
\theoremsymbol{\rule{1ex}{1ex}}

\theoremstyle{empty} \theoremseparator{.}

\def\cal{\mathcal}

\def\Bbb{\mathbb}
\def\bar{\overline}

\def\N{\Bbb{N}}
\def\Naturals{\Bbb{N}}
\def\Reals{\Bbb{R}}

\def\ni{\noindent}

\def\Area{\hbox{\rm Area}}

\def\Time{\hbox{\rm Time}}
\def\Space{\hbox{\rm Space}}

\def\Height{\hbox{\rm Height}}

\def\FL{\hbox{\rm FL}}

\def\ssm{\smallsetminus}

\def\ms{\medskip}

\def\onto{{\kern3pt\to\kern-8pt\to\kern3pt}}

\def\<{\langle}
\def\>{\rangle}
\def\|{{\ |\ }}

 \def\inv{^{-1}}

 \def\AA{\cal A}

 \def\RR{\cal R}

 \def\PP{\cal P}
\def\QQ{\cal Q}
\def\XX{\cal X}
\def\YY{\cal Y}
 
\def\inv{^{-1}}

\newcommand{\aka}{\emph{a.k.a.}\ }
\newcommand{\set}[1]{\left\{#1\right\}}

\newcommand{\abs}[1]{\left|#1\right|}
\renewcommand{\ni}{\noindent}

\renewcommand{\ms}{\medskip}
\newcommand{\bs}{\bigskip}

\def\qed{{\ifhmode\unskip\nobreak\hfil\penalty50 \hskip1em \else\nobreak\fi
   \mbox{}\nobreak\hfil \rule{1ex}{1ex}%
   \parfillskip=0pt \finalhyphendemerits=0 \par}}

\begin{document}

\title{Isoperimetric inequalities for nilpotent groups}

\author{S. M. Gersten,  D. F. Holt
 and T. R. Riley}
\date{January 16, 2002}
\maketitle
\begin{abstract}

We prove that every finitely generated nilpotent group of class
$c$ admits a polynomial isoperimetric function of degree $c+1$ and
a linear upper bound on its filling length function.

\medskip
\footnotesize{\ni \textbf{1991 Mathematics Subject
Classification:} 20F05, 20F32, 57M07 \\ \ni \emph{Key words and
phrases:} nilpotent group, isoperimetric function, filling length}
\end{abstract}

\section{Introduction}

The main result of this article is a proof of what has been known
as the \emph{``$c+1$-conjecture''} (see \cite{Bridson_survey} or
$5\textup{A}'_5$ of \cite{Gromov}).
\begin{MainThm}[$c+1$-conjecture]\label{c+1-theorem}
Every finitely generated nilpotent group $G$ admits a polynomial
isoperimetric inequality of degree $c+1$, where $c$ is the
nilpotency class of $G$.
\end{MainThm}
An isoperimetric inequality for a finite presentation $\mathcal{P}
= \langle \AA \mid  \RR \rangle$ of a group $G$ concerns
\emph{null-homotopic} words $w$, that is words that evaluate to
the identity in $G$. It gives an upper bound (an
\emph{``isoperimetric function''}), in terms of the length of $w$,
on the number of times one has to \emph{apply relators} from $\RR$
to $w$ in a process of reducing it to the empty word. (More
details are given in \S2.) The following couple of remarks are
important in making sense of the statement of
Theorem~\ref{c+1-theorem}. Nilpotent groups are \emph{coherent},
by which we mean that their finitely generated subgroups are
finitely presentable.  In particular, all finitely generated
nilpotent groups have a finite presentation. Further (see \S2) an
isoperimetric function $f_{\mathcal{P}}(n)$ concerns a fixed
finite presentation $\mathcal{P}$ for $G$; however if
$\mathcal{Q}$ is another finite presentation for $G$ then there is
an isoperimetric function $f_{\mathcal{Q}}(n)$ for $\QQ$ that is
that satisfies $f_{\mathcal{P}} \simeq f_{\mathcal{Q}}$, where
$\simeq$ is a well-known equivalence relation (see
Definition~\ref{simeq}). Moreover, if $f_{\mathcal{P}}$ is
polynomial of some given degree $\geq 1$, then we can always take
$f_{\mathcal{Q}}$ to be a polynomial of the same degree.

\ms Our strategy will be to prove the $c+1$-conjecture by an
induction on the class $c$. However we use an induction argument
in which we keep track of more than an isoperimetric function. In
fact, we prove the following stronger theorem.

\begin{MainThm}\label{AFL-pair thm} Finitely generated nilpotent groups
$G$ of class $c$ admit $(n^{c+1},n)$ as an $(\Area, \FL)$-pair.
\end{MainThm}

The terminology used above is defined carefully in \S2, but
essentially, what this theorem says is the following.  Suppose $w$
is a null-homotopic word $w$ of length $n$ in some finite
presentation for $G$. We can reduce $w$ to the empty word by
applying at most $O(n^{c+1})$ relators from the presentation, and
in such a way that in the process the intermediate words have
length at most $O(n)$.

Let $$G \ = \ \Gamma_1 > \Gamma_2 > \cdots > \Gamma_{c+1} =
\set{1}$$ be the lower central series for $G$ defined inductively
by $\Gamma_1 := G$ and $\Gamma_{i+1} := [ G, \Gamma_i]$.

We will use a generating set $\AA$ for $G$ which will be a
disjoint union of sets $\AA_i$, where, for each $i$, $\AA_i$ is a
generating set for $\Gamma_i$ modulo $\Gamma_{i+1}$.  For $\AA_1$,
we take inverse images in $G$ of generators of cyclic invariant
factors of the abelian group $G / \Gamma_2$.  Then, inductively,
for $i > 0$, we define $\AA_{i+1} := \set{ [x,y] \mid x \in \AA_1,
y \in \AA_i }$.

The idea of the proof of Theorem~\ref{AFL-pair thm} is to start by
reducing the word $w$ to the identity modulo $\Gamma_c$. The
quotient $G/ \Gamma_c$ is nilpotent of class $c$ and so, by
induction hypothesis, we can reduce $w$ to the identity in a
presentation for $G/ \Gamma_c$, by applying $O(n^c)$ defining
relators. When we carry out the corresponding reduction in a
presentation for $G$ we introduce $O(n^c)$ extra generators from
$\AA_c$. We use the definition of these generators as commutator
 words $z$ to compress their powers $z^m$ with $m = O(n^c)$ to words of
length $O(n)$. This compression process is handled in
Proposition~\ref{compressing commutators} and
Corollary~\ref{area}, and we show that it can be accomplished by
applying $O(n^{c+1})$ relators.

\ms In addition to proving the $c+1$-conjecture,
Theorem~\ref{AFL-pair thm} yields the following corollary.

\begin{MainCor}\label{fl cor} If $\PP$ is a finite presentation
for a nilpotent group then there exists $\lambda >0$ such that the
filling length function $\FL$ of $\PP$ satisfies $\FL(n) \leq
\lambda n$ for all $n \in \Naturals$.
\end{MainCor}

This result was proved by the third author in \cite{Riley} via the
rather indirect technique of using asymptotic cones: finitely
generated nilpotent groups have simply connected asymptotic cones
(see Pansu \cite{Pansu}) and groups with simply connected
asymptotic cones have linearly bounded filling length.  This
article provides a direct combinatorial proof.  Note, also, that
the filling length function for a finite presentation $\PP$ is an
isodiametric function for $\PP$.  So we also have a direct
combinatorial proof that finitely generated nilpotent groups admit
linear isodiametric functions.

\ms

There is a computer science reinterpretation of
Theorem~\ref{AFL-pair thm} because $\Area$ and $\FL$ can be
recognised to be measures of computational complexity in the
following context.

Suppose $\PP$ is a finite presentation for a group $G$.  One can
attempt to solve the word problem using a non-deterministic Turing
machine as follows. Initially the input word $w$ of length $n$ is
displayed on the Turing tape. A \emph{step} in the operation of
the machine is an application of a relator, a free reduction, or a
free expansion (see \S2).  The machine searches for a \emph{proof}
that $w=1$ in $G$ -- that is, a sequence of steps that reduces $w$
to the empty word.  The running time of a \emph{proof} is the
number of steps, and its space is the number of different entry
squares on the Turing tape that are disturbed in the course of the
\emph{proof}.

The time $\Time(w)$ for a word $w$ such that $w = 1$ in $G$ is the
minimum running time amongst \emph{proofs} for $w$, and the time
function $\Time: \Naturals \to \Naturals$ is defined by $\Time(n)
:= \max \set{ \Time(w) \mid \ell(w) \leq n \text{ and } w=1 \text{
in } G}$. Similarly we define the space function $\Space$ by
setting $\Space(w)$ to be the minimal space amongst all proofs for
$w$.

On comparing with the definitions of \S2, we quickly recognise
that the space function $\Space$ is precisely the same as the
filling function $\FL$. Also the time function $\Time$ is closely
related to the function $\Area$. (In fact $\Time$ is precisely the
function $\Height$ of Remark~\ref{Area-Height}, where it is
explained how $\Height$ relates to $\Area$.)  Thus we have the
following corollary to Theorem~\ref{AFL-pair thm}.

\begin{MainCor}
Given a finite presentation for a nilpotent group $G$ of class
$c$, the non-deterministic Turing machine described above solves
the word problem with $\Space(n)$ bounded above by a linear
function of $n$ and $\Time(n)$ bounded above by a polynomial in
$n$ of degree $c+1$. Moreover, given an input word $w$ such that
$w=1$ in $G$, there is a \emph{proof} for $w$ that runs within
both these bounds simultaneously.
\end{MainCor}

\noindent
[Note that these complexity bounds are not the best possible.
The torsion subgroup $T$ of a finitely generated nilpotent
group $G$ is finite, and by a result of S.A.~Jennings
(see~\cite{Hall} or Theorem~2.5 of~\cite{Wehrfritz}),
the torsion-free group $G/T$ can be embedded in an upper
unitriangular matrix group $U$ over $\mathbb{Z}$.
It is easy to see that the integer entries of a matrix in $U$ representing
a word of length $n$ in $G/T$ are at most $O(n^{d-1})$, where $d$ is
the degree of the matrix. Using this and the fact that
two $k$-digit numbers can be multiplied in time
$O(k\, (\log k)^2 )$, we can obtain a deterministic solution to the word
problem in $G$  in time  $O(n \, (\log n)^2 )$ and space $O ( \log n)$.

Even this is not the best possible result on the time complexity.
In the unpublished manuscript~\cite{CGS}, Cannon, Goodman and Shapiro
show that finitely generated nilpotent groups admit a Dehn algorithm if
we adjoin some extra symbols to the group generators, and this leads
to a linear solution of the word problem in $G$.]

\ms The proof of the $c+1$-conjecture is the culmination of a
number of results spanning the last decade or so. The first author
proved in \cite{Gersten6} that $G$ admits a polynomial
isoperimetric inequality of degree $2^h$, where $h$ is the Hirsch
length of $G$. The degree was improved to $2 \cdot 3^c$ by
G.~Conner \cite{Conner}, and then improved further to $2c$ by
C.~Hidber \cite{Hidber}.

We also mention Ch.~Pittet \cite{Pittet}, who proved that a
lattice in a simply connected homogeneous Lie group of class $c$
admits a polynomial isoperimetric function of degree $c+1$. (The
nilpotent Lie group is called \emph{homogeneous} if its Lie
algebra is graded.) Gromov also suggested a possible reduction to
the homogeneous case by perturbing the structure constants
\cite{Gromov} $5.\textup{A}_5$, but no-one has yet succeeded in
carrying out this plan.

The isoperimetric inequality proved in this article is the best
possible bound in terms of the class in general.  For example if
$G$ is a free nilpotent group of class $c$ then its minimal
isoperimetric function (\aka its Dehn function) is polynomial of
degree $c+1$ (see \cite{BMS} or \cite{Gersten} for the lower bound
and \cite{Pittet} for the upper bound). However it is not best
possible for individual nilpotent groups: for example D.~Allcock
in \cite{Allcock} gives the first proof that the $2n+1$
dimensional integral Heisenberg groups admit quadratic
isoperimetric functions for $n>1$; these groups are all nilpotent
of class 2.  By way of contrast, the 3-dimensional integral
Heisenberg group has a cubic minimal isoperimetric polynomial
(\cite{Epstein} and \cite{Gersten4}).

\medskip
We give two corollaries of Theorem~\ref{c+1-theorem} concerning
the cohomology of groups and
 differential geometry
respectively.

The first corollary is about the growth of cohomology classes. Let
$G$ be a finitely generated group with word metric determined by a
finite set of generators and denoted $d(1,g) = \abs{g}$. Recall
that a real valued 2-cocycle (for the trivial $G$-action on
$\mathbb{R}$) on the bar construction is a function $f : G \times
G \to \Reals$ satisfying $f(x,y z)+f(y,z)=f(x,y)+f(x y,z)$ for all
$x,y,z \in G$.

\addtocounter{MainThm}{-1} %
\addtocounter{MainCor}{-2}

\begin{MainCor}
Let $G$ be a finitely generated nilpotent group of class $c$ and
let $\zeta \in H^2(G, \Reals)$.  Then there is a 2-cocycle $f$ in
the class $\zeta$ which satisfies $f(x,y) \leq
M(\abs{x}+\abs{y})^{c+1}$ for all $x,y \in G$ and constant $M >0$.
\end{MainCor}

The terminology (due to Gromov) of an \emph{isoperimetric
function} for a group is motivated by the analogous notion with
the same name from differential geometry.  One can draw parallels
between van~Kampen diagrams (see Remark~\ref{different defns})
filling edge-circuits in the Cayley 2-complex associated to a
finite presentation of a group and homotopy discs for rectifiable
loops in the universal cover of a Riemannian manifold.

\begin{MainCor} \label{diff geom cor}
Let $M$ be a closed Riemannian manifold whose fundamental group is
nilpotent of class $c$.  Then there is a polynomial $f$ of degree
$c+1$ such that for every rectifiable loop of length $L$ in the
universal cover $\tilde{M}$, there is a singular disc filling of
area at most $f(L)$.
\end{MainCor}

For the proof one considers a piecewise $C_1$-map of the
presentation complex $K$ of a finite presentation of $G$ into $M$
and the induced map of universal covers $\tilde{K} \to \tilde{M}$.
An edge-circuit $p$ in $\tilde{K}$ defines a null-homotopic word
$w$ in $\pi_1 M$, and a van~Kampen diagram $D$ for $w$ gives a
Lipschitz filling for $p$ in $\tilde{K}$.   For a general
Lipschitz loop $p$ in $\tilde{M}$, one uses Theorem~10.3.3 of
\cite{Epstein} to homotop $p$ by a Lipschitz homotopy to an
edge-circuit in $\tilde{K}$, which is then filling to realise the
desired isoperimetric inequality for $p$.


\section{Isoperimetric functions and filling
 length functions}\label{Dehn defn}

Let $\PP=\langle \AA \mid \RR \rangle$ be a finite presentation
for a group $G$.  A \emph{word} $w$ is an element of the free
monoid $(\mathcal{A}\cup\mathcal{A}^{-1})^{\star}$.  Denote the
length of $w$ by $\ell(w)$. We say $w$ is \emph{null-homotopic}
when $w=1$ in $G$.  For a word $w = {a_1}^{\varepsilon_1}
{a_{2}}^{\varepsilon_2} \dots {a_s}^{\varepsilon_s}$, where each
$a_i \in \AA$ and each $\varepsilon_i = \pm 1$, the \emph{inverse
word} $w^{-1}$ is ${a_s}^{-\varepsilon_s}
 \dots {a_2}^{-\varepsilon_2} {a_1}^{-\varepsilon_1} $.

\begin{defn}  A \emph{$\PP$-sequence} $S$ is a finite sequence of words $w_0,
w_1,\dots, w_m$ such that each $w_{i+1}$ is obtained from $w_i$ by
one of three moves.
\begin{enumerate}
\item \emph{Free reduction.} Remove a subword $a
a^{-1}$ from $w_i$, where $a$ is a generator or the inverse of a
generator.
\item \emph{Free expansion.} Insert a subword $a
a^{-1}$ into $w_i$, where $a$ is a generator or the inverse of a
generator. So $w_{i+1}=uaa^{-1}v$ for some words $u, \, v$ such
that $w_i=uv$ in  $(\mathcal{A}\cup\mathcal{A}^{-1})^{\star}$.
\item \emph{Application of a relator.} Replace $w_i=\alpha u\beta$ by
$\alpha v\beta$, where $uv\inv$ is a cyclic conjugate of one of
the defining relators or its inverse.
\end{enumerate}
We refer to $m$ as the \emph{height} $\Height(S)$ of $S$ and we
define the \emph{filling length} $\FL(S)$ of $S$ by $$\FL(S):=
\max \set{\ell(w_i) \mid 0 \leq i \leq m}.$$  The \emph{area}
$\Area(S)$ of $S$ is defined to be the number of $i$ such that
$w_{i+1}$ is obtained from $w_i$ by an \emph{application of a
relator} move.

If $w=w_0$ and $w_m$ is the empty word then we say that $S$ is a
\emph{null-$\PP$-sequence} for $w$ (or, more briefly, just a
``\emph{null-sequence}'' when $\PP$ is clear from the context).
\end{defn}

\ms In this article we are concerned with two \emph{filling
functions}, that measure ``\emph{area}'' and ``\emph{filling
length}'' -- two different aspects of the \emph{geometry of the
word problem} for $G$.

\begin{defn}  Let $w$ be a null-homotopic word in $\PP$.

We define the \emph{area} $\Area(w)$ of $w$ by $$\Area(w):=
\min\set{\Area(S) \mid \textup{null-}\PP\textup{-sequences } S
\text{ for } w },$$ that is, the minimum number of relators that
one has to apply to reduce $w$ to the empty word.  Similarly, the
\emph{filling length} of $w$ is $$\FL(w):= \min \set{\FL(S) \mid
\textup{null-}\PP\textup{-sequences } S \textup{ for } w }.$$

We define the \emph{Dehn function} $\Area:\Naturals \to \Naturals$
(also known as the \emph{minimal isoperimetric function}) and the
\emph{filling length function} $\FL:\Naturals \to \Naturals$ by
\begin{eqnarray*}
\Area(n) & := & \max \set{\Area(w) \mid \textup{null-homotopic
words } w \textup{ with } \ell(w) \leq n } \\ \FL(n) & := & \max
\set{\FL(w) \mid \textup{null-homotopic words } w \textup{ with }
\ell(w) \leq n }.
\end{eqnarray*}
\end{defn}

\begin{Remark} \label{different defns}
The formulations of the definitions above are those we will use in
this article, but we mention some equivalent alternatives that
occur elsewhere in the literature.

We could equivalently define $\Area(w)$ to be the minimal number
of 2-cells needed to construct a van~Kampen diagram for $w$, or
the minimal $N$ such that there is an equality $w= \prod_{i=1}^{N}
u_i^{-1} r_i u_i$ in the free group $F(\AA)$ for some $r_i \in
\RR^{\pm 1}$ and words $u_i$.

Similarly, $\FL(w)$ can be defined in terms of \emph{shellings} of
van~Kampen diagrams $D_w$ for $w$.  A \emph{shelling} of $D_w$ is
a combinatorial null-homotopy down to the base vertex.  The
filling length of a shelling is the maximum length of the boundary
loops of the diagrams one encounters in the course of the
null-homotopy.  The filling length of $D_w$ is defined to be the
minimal filling length amongst shellings of $D_w$ and then
$\FL(w)$ can be defined to be $$\min\set{\FL(D_w) \mid
\textup{van~Kampen diagrams } D_w \text{ for } w}.$$  For detailed
definitions and proofs we refer the reader to
\cite{Gersten-Riley}.
\end{Remark}

The following is a technical lemma that we will use in \S3.

\begin{lemma} \label{app of a relator}
Let $S$ be a $\PP$-sequence $w_0, w_1, \ldots, w_m$ as defined
above. Let $C:= \max \set{\ell(r) \mid r \in \RR}$.  There is
another $\PP$-sequence $\hat{w}_0, \hat{w}_1, \ldots,
\hat{w}_{\hat{m}}$, which we will call $\hat S$, such that
\begin{itemize}
\item  $\hat {w}_0 = w_0$, $\hat {w}_{\hat{m}} = w_m$
\item $\Area(\hat{S}) = \Area(S) $
\item $\FL(\hat{S}) \leq \FL(S) + C $
\end{itemize}
and such that every time $\hat{w}_{i+1}$ is obtained from a word
$\hat{w}_i$ in the sequence $\hat{S}$ by an \emph{application of a relator}
the whole of a cyclic conjugate of an element of $\RR$ or its
inverse is inserted into $\hat{w}_i$.
\end{lemma}

\ni \emph{Proof.} Suppose that $w_{i+1}$ is obtained from a word
$w_i$ in the sequence $S$ by an \emph{application of a relator}:
that is, $w_i=\alpha u \beta$ and $w_{i+1}=\alpha v \beta$ for
some words $\alpha, \beta, u,v $ where $u v\inv$ is a cyclic
conjugate some element of $\RR^{\pm 1}$. We can obtain $w_{i+1}$
from $w_i$ by inserting ${u}^{-1}v$ into $w_i= \alpha u \beta$ to
get $\alpha {u}{u}^{-1}{v}\beta$ and then using at most $C$ free
reductions to retrieve $\alpha {v} \beta$. \qed

\ms

\begin{Remark}\label{Area-Height}
We mention that the area of a null-homotopic word $w$ is closely
related to the \emph{height} of its null-sequences. Define
$$\Height(w) \ := \ \min \set{\Height(S) \mid
\textup{null-sequences } S \text{ for } w }.$$ Then for all
null-homotopic words $w$, $$\Area(w) \ \leq \ \Height(w) \ \leq \
(C+1) \Area(w)+ \ell(w),$$ where $C  := \max \set{ \ell(r) \mid r
\in \RR }$.

The inequality $\Area(w) \leq \Height(w)$ follows from the
definitions since free reductions and expansions do not contribute
to the area. To obtain the inequality $\Height(w) \leq (C+1)
\Area(w)+ \ell(w)$ first take a van~Kampen diagram $D_w$ for $w$
with $\Area(D_w)= \Area(n)$. Take any shelling $$D_w \ = \ D_0,
D_1, \ldots, D_m  \ = \ \star$$ of $D_w$ down to its base vertex
$\star$ in which each $D_{i+1}$ is obtained from $D_i$ by a
\emph{2-cell collapse} or a \emph{1-cell collapse} (but never a
\emph{1-cell expansion}) and let $w_j$ be the boundary word of
$D_j$.  Then $w_0,w_1, \ldots, w_m$ is a null-sequence for $w$,
where $w_{i+1}$ is obtained from $w_i$ by applying of a relator if
$D_{i+1}$ is obtained from $D_i$ by a \emph{2-cell collapse}, and
by a free reduction if $D_{i+1}$ is obtained from $D_i$ by a
\emph{1-cell collapse}. We obtain the required inequality by
observing that the number of \emph{2-cell collapse} moves in
the shelling is $\Area(w)$ and the total number of \emph{1-cell
collapse} moves is at most the total number if 1-cells in the
1-skeleton of $D_w$, which is at most $C \Area(w)+ \ell(w)$.
\end{Remark}

An \emph{isoperimetric inequality} for $\PP$ is provided by any
function $f : \N \to \N$ such that $\Area(n) \leq f(n)$ for all
$n$.  We refer to $f$ as an \emph{isoperimetric function} for
$\PP$.

\medskip
It is important to note that $\Area: \Naturals \to \Naturals$ and
$\FL: \Naturals \to \Naturals$ are both defined for a fixed finite
presentation.  However each function is a group invariant in the
sense that each is well behaved under change of finite
presentation as we now explain.
\begin{defn} \label{simeq}
For two functions $f,g : \N \to \N$ we say that $f \preceq g$ when
there exists $C>0$ such that $f(n) \leq C g(Cn+C) + C n +C$ for
all $n$, and we say $f \simeq g$ if and only if $f \preceq g$ and
$g \preceq f$.
\end{defn}
Let $\PP$ and $\QQ$ be two presentations of the same group $G$,
let $\Area_{\PP}$ and $\Area_{\QQ}$ be their Dehn functions, and
let $\FL_{\PP}$ and $\FL_{\QQ}$ be their filling length functions.
Then it is proved in \cite{Gersten-Short} that $\Area_{\PP} \simeq
\Area_{\QQ}$ and in \cite{Gersten-Riley} that $\FL_{\PP} \simeq
\FL_{\QQ}$.

\begin{defn}
We say that a pair $(f,g)$ of functions $\Naturals \to \Naturals$
is an $(\Area, \FL)$-pair for $\PP$ when there exists a constant
$\lambda >0$ such that for any null-homotopic word $w$, there
exists a null-$\PP$-sequence $S$ with
\begin{eqnarray*}
\Area(S) & \leq & \lambda f(\ell(w)) \\ \FL(S) & \leq & \lambda
g(\ell(w)).
\end{eqnarray*}
\end{defn}

The first part of the following lemma justifies the deduction of
Theorem~\ref{c+1-theorem} and Corollary~\ref{fl cor} from
Theorem~\ref{AFL-pair thm}, and second part  justifies the
terminology used in the statement of Theorem~\ref{AFL-pair thm}.

\begin{lemma}\label{AFL-pair lemma}
If $(f,g)$ is an $(\Area, \FL)$-pair for $\PP$ then  $f$ is an
isoperimetric function for $\PP$ and $g$ is an upper bound for its
filling length function.

Moreover, if $(n^\alpha,n^\beta)$ is an $(\Area, \FL)$-pair
for
$\PP$ where $\alpha,\beta\ge 1$, then
$(n^\alpha,n^\beta)$ is also an $(\Area, \FL)$-pair for
any other finite presentation $\QQ$ of the same
group $G$.

\end{lemma}

\ni \emph{Proof.} The first part of the lemma follows immediately
from the definitions.

From the proofs given in \cite{Gersten-Riley}and
\cite{Gersten-Short}
 it follows that if
$(f,g)$ is an $(\Area, \FL)$-pair for $\PP$, then $(f',g')$ is an
$(\Area, \FL)$-pair for $\QQ$, where $f'(n)=Cf(Cn)+Cn$ and
$g'(n)=Cg(Cn)+Cn$ for some $C>0$; here we have eliminated the
additive constants in the definition of the equivalence relation
$\simeq$ by making use of the fact that the empty word is the only
one of length~0.

But if $f(n)=n^\alpha$ with $\alpha\ge 1$,
then $Cf(Cn)+Cn\le C'n^\alpha$
for suitable constant $C'>0$.
It follows that if $(n^\alpha,n^\beta)$ is
an $(\Area, \FL)$-pair for $\PP$, then
it is also an $(\Area, \FL)$-pair for $\QQ$.
\qed

\section{Proof of Theorem~\ref{AFL-pair thm}}

\ms We say that an element of a nilpotent group $G$ has
\emph{weight} $k$ if it lies in $\Gamma_k \smallsetminus
\Gamma_{k+1}$, where $\{ \Gamma_i \}$ is the lower central series
of $G$, as defined in \S1. Before we come to the proof of
Theorem~\ref{AFL-pair thm} we give a proposition and corollary
concerning \emph{compressing} powers of elements of weight $c$ in
a finite presentation for a nilpotent group of class $c$.

We use the following conventions for commutator words:
\begin{eqnarray*}
[a] & := & a \\ %
\mbox{$[ \, a \, , \, b \, ]$} & := &  a^{-1} \, b^{-1} \, a  \, b \\ %
\mbox{$[ \, a_1 \, , \, a_2 \, , \, \ldots, \, a_{c-1} \, , \, a_c
\,  ]$} & := & [ \, a_1 \, , \, [ \, a_2 \, , \, \ldots , [ \,
a_{c-1} \, , \, a_c \, ] \ldots ]].
\end{eqnarray*}

Let  $\XX = \set{ x_1, x_2, \ldots, x_c }$ be an alphabet and for
$k=1, 2, \ldots, c$ define
\begin{eqnarray*}
\XX_k & := & \set{x_k, x_{k+1}, \ldots, x_c} \\ \RR_k & := &
\set{[y_k, y_{k+1}, \ldots, y_{c+1}] \mid y_j \in {\XX_k}^{\pm 1},
\, k \le j \le c+1  \, }.
\end{eqnarray*}
Let $\PP_k := \langle \XX_k \mid \RR_k \rangle$, which is a finite
presentation for a free nilpotent group $G_k$ of class $c +1 - k
$.  Define $z_k$ to be the word $[ \, x_k \, , \, x_{k+1} \, , \,
\ldots , \, x_c \, ]$. In particular, we have $z_c = x_c$.

One can regard the first part of the following proposition as
comparing the use of commutators in a free nilpotent group with
the representation of positive integers $s$ by their $n$-ary
expansion $s= s_0 + s_1 n + s_2 n^2 + \ldots$.   The second part
gives upper bounds on the area and filling length of a
$\PP$-sequence between the commutators corresponding to the
$n$-ary representations for $s$ and $s+1$.

\begin{prop}\label{compressing commutators}
Suppose $c \in \Naturals \ssm \set{0}$ and $s,n \in \Naturals$
with $0 \leq s \leq n^c-1$. Express $s$ as a sum: $$  s \ = \ s_0
\ + \ s_1 \, n \ + \ \ldots \ + \ s_{c-1} \, n^{c-1},$$ where each
$s_i \in \set{0, 1, \ldots ,n-1}$.  Then, with the notation above,
the word $$\widetilde{{z_1}^s} \ := \ {z_1}^{s_0} \, [ \, {x_1}^n
\, , \,{z_2}^{s_1} \,[{x_2}^n  \, , \, \ldots , \,
{z_{c-1}}^{s_{c-2}} \, [ \, {x _{c-1}}^n \, , \, {z_c}^{s_{c-1}}
\, ] \ldots ] \, ]$$ equals ${z_1}^s$ in $G_1$. Also
$$\widetilde{{z_1}^{n^c}} \ := \ [ \, {x_1}^n \, , \, {x_2}^n \, ,
\, \ldots , \, {x_c}^n \, ]$$ equals ${z_1}^{n^c}$ in $G_1$.

Moreover, there exists $\kappa>0$ such that for $0 \leq s \leq
n^c-1$ we can transform $z_1 \widetilde{{z_1}^s}$ to
$\widetilde{{z_1}^{s+1}}$ via a $\PP_1$-sequence of filling length
at most $\kappa n$ and of area at most $\kappa n^{k+1}$ where $k$
is the integer such that $n^ k \mid (s+1)$ but $n^{k+1} \nmid
(s+1)$.
\end{prop}

\ms \ni \emph{Proof.} Our proof is by induction on the integer $c$.
The base case of $c=1$ is straightforward. We just take
$\widetilde{{z_1}^s}:= {z_1}^s = {x_1}^s$ and $\kappa:= 1$.

We now prove the induction step.  Suppose $ c \in \Naturals \ssm
\set{0,1}$. Since the free reduction of $\widetilde{{z_1}^0}$
is the empty word, we have $\widetilde{{z_1}^0} =1$ in $G$.
It therefore suffices to show that for $s,n \in \Naturals$
with $1 \leq s \leq n^c-1$, there is a $\PP_1$-sequence from $z_1
\widetilde{{z_1}^s}$ to $\widetilde{{z_1}^{s+1}}$ within the
required filling length and area bounds.

Express $s$ as a sum: $$  s \ = \ s_0 \ + \ s_1 \, n \ + \ \ldots
\ + \ s_{c-1} \, n^{c-1},$$ where each $s_i \in \set{0, 1, \ldots
,n-1}$.  Define $$t \ := \ s_1 \ + \ s_2 \, n \ + \ \ldots \ + \
s_{c-2} \, n^{c-2}.$$  Then
$\widetilde{{z_1}^s} \ = \ {z_1}^{s_0} \, [ \, {x_1}^n \, , \,
\widetilde{{z_2}^t} \, ]$ as words.

If $s_0+1<n$ then $z_1 \widetilde{{z_1}^s} =
\widetilde{{z_1}^{s+1}}$ as words and there is a trivial
$\PP_1$-sequence from $z_1 \widetilde{{z_1}^s}$ to
$\widetilde{{z_1}^{s+1}}$. If $s_0+1=n$ then we calculate that
$z_1 \widetilde{{z_1}^s} = \widetilde{{z_1}^{s+1}}$ in $G_1$ as
follows:
\begin{eqnarray}
z_1 \widetilde{{z_1}^s} & = & {z_1}^n \, [ \, {x_1}^n \, ,
\, \widetilde{{z_2}^t} \, ] \\ %
& = & {z_1}^n \, {x_1}^{-n} \, {(\widetilde{{z_2}^t})}^{-1} \,
{x_1}^n \, \widetilde{{z_2}^t} \\ %
& = & {z_1}^n \, {x_1}^{-n} \, {(\widetilde{{z_2}^t})}^{-1} \,
{z_2}^{-1} \, z_2 \, {x_1}^n   \, \widetilde{{z_2}^t} \\ %
& = & {x_1}^{-n} \, {(\widetilde{{z_2}^t})}^{-1} \, {z_2}^{-1} \,
{x_1}^n \, z_2 \, \widetilde{{z_2}^t} \\ %
& = & [ \, {x_1}^n \, , \, {z_2} \, \widetilde{{z_2}^t} \, ] \\ %
& = & [ \, {x_1}^n \, , \, \widetilde{{z_2}^{t+1}} \, ] \\ %
& = & \widetilde{{z_1}^{s+1}}.
\end{eqnarray}
In the step from (3) to (4) we use the fact that $z_1  \ = \ [ \,
x_1 \, , \, z_2 \, ]$ and is central.   For the step from (5) to
(6) we invoke the induction hypothesis to tell us that ${z_2} \,
\widetilde{{z_2}^t} \ = \ \widetilde{{z_2}^{t+1}} $ in $G_2$, from
which it follows that ${z_2} \, \widetilde{{z_2}^t} \,
{(\widetilde{{z_2}^{t+1}})}^{-1}$ is central in $G_1$ and $[ \,
{x_1}^n \, , \, {z_2} \, \widetilde{{z_2}^t} \, ] \ = \ [ \,
{x_1}^n \, , \, \widetilde{{z_2}^{t+1}} \, ]$.

The course of the calculation above dictates how to construct a
$\PP_1$-sequence from $z_1 \widetilde{{z_1}^s}$ to
$\widetilde{{z_1}^{s+1}}$. The steps that require some further
explanation are those at which \emph{application of relator} moves
are used: that is, from (3) to (4) and from (5) to (6).

The step from (3) to (4) is performed by a $\PP_1$-sequence in
which we introduce $n$ commutator words $ {z_1}^{-1} = [x_1 ,
z_2]^{-1}$ to move $z_2$ past ${x_1}^n$.  After a word
${z_1}^{-1}$ is introduced it is immediately moved through the
word using relations from $\RR_1$ and is then cancelled with a
$z_1$. The number of letters each $z_1$ has to be moved past is
bounded by $n$ up to a multiplicative constant that depends only
on $\PP_1$ (since $\ell(\widetilde{{z_2}^t}) \leq \bar{\kappa} n$,
where $\bar{\kappa}>0$ is a constant depending only on $\PP_2$).
So the area of this $\PP_1$-sequence is at most $n^2$ and its
filling length is at most $n$ up to a multiplicative constant.

Now let us explain how to construct a $\PP_1$-sequence for the
step from (5) to (6).  Suppose that $n^k \mid (s+1)$ but $n^{k+1}
\nmid (s+1)$.  Then $n^{k-1} \mid (t+1)$ but $n^{k} \nmid (t+1)$.
By induction hypothesis there is a $\PP_2$-sequence $\bar S $ from
$z_2 \, \widetilde{{z_2}^t}$ to $\widetilde{{z_2}^{t+1}}$ with
area at most $\bar{\kappa} n^{k}$ and filling length at most
$\bar{\kappa} n$.

Lemma~\ref{app of a relator} allows us to assume (by suitably
adjusting the constant $\bar{\kappa}$) that in every instance of
an \emph{application of a relator} move in the sequence $\bar{S}$ a
whole cyclic conjugate of a relator is inserted.

We now explain how to use $\bar{S}$ to induce a $\PP_1$-sequence
that transforms $[ \, {x_1}^n \, , \, {z_2} \, \widetilde{{z_2}^t}
\, ] $ into $[ \, {x_1}^n \, , \, \widetilde{{z_2}^{t+1}} \, ]$.
First define $\bar{S}'$ to be the $\PP_2$-sequence that transforms
${(\widetilde{{z_2}^t})}^{-1} \, {z_2}^{-1}$ to
${(\widetilde{{z_2}^{t+1}})}^{-1}$ and is obtained by inverting
every word in $\bar{S}$.  Now $$[ \, {x_1}^n \, , \, {z_2} \,
\widetilde{{z_2}^t} \, ] = {x_1}^{-n} \,
{(\widetilde{{z_2}^t})}^{-1} \, {z_2}^{-1} \, {x_1}^n \,  {z_2} \,
\widetilde{{z_2}^t}.$$ Consider running the $\PP_2$-sequences
$\bar{S}$ and $\bar{S}'$ concurrently on the subwords ${z_2} \,
\widetilde{{z_2}^t}$ and ${(\widetilde{{z_2}^t})}^{-1} \,
{z_2}^{-1}$ respectively in $[ \, {x_1}^n \, , \, {z_2} \,
\widetilde{{z_2}^t} \, ]$: that is, we do the first move in
$\bar{S}$ and then the first move in $\bar{S}'$, then the second
move in $\bar{S}$ and then the second move in $\bar{S}'$, and so
on. However we want to construct a $\PP_1$-sequence, not a
$\PP_2$-sequence.

Suppose a move in $\bar{S}$ is the insertion of a word $r \in
\RR_2$.  Then the corresponding move in $\bar{S}'$ inserts the
word $r^{-1}$. Use free expansion moves to insert the word
$r^{-1}r$ in the place where $\bar{S}$ dictated that $r$ was to be
inserted; then use  relators in $\RR_1$ to move $r^{-1}$ through
the word to the place where $\bar{S}'$ dictated $r^{-1}$ was to be
inserted.  (Recall that $r$ represents a central element of $G_1$
and the appropriate commutator to move it past letters from
$\XX_1$ are in $\RR_1$.)

The number of relators from $\RR_1$ that have to be applied to
move each $r^{-1}$ to its appropriate place is at most $n$ up to a
multiplicative constant. It follows that we can find a constant
$\kappa>0$ such that there is a $\PP_1$-sequence that transforms
the word $z_1 \widetilde{{z_1}^s}$ into $\widetilde{{z_1}^{s+1}}$,
with area at most $\kappa n^{k+1}$ and filling length at most
$\kappa n$.  This completes the proof. \qed

\bs

Recall from \S1 that we define a set of generators $\AA$ for a
 finitely generated nilpotent group $G$ of class $c$, as follows.
The set $\AA$ is a disjoint union of sets $\AA_i$, where, for each
$i$, $\AA_i$ is a generating set for $\Gamma_i$ modulo
$\Gamma_{i+1}$. For $\AA_1$, we take inverse images in $G$ of
generators of cyclic invariant factors of the abelian group $G /
\Gamma_2$.  Then, inductively, for $i > 0$, we define $\AA_{i+1}
:= \set{ [x,y] \mid x \in \AA_1, y \in \AA_i }$.

\begin{cor}\label{area}
Let $G$ be a finitely generated nilpotent group of class $c$, and
let $\AA$ be the generating set defined above.  Then, for any
finite presentation $\PP = \langle \AA \mid \RR \rangle$ for $G$,
there is a constant $\xi>0$ depending only on $\PP$ with the
following properties.

Let $n \in \Naturals$. For each $z \in \AA_c$ there are
``\emph{compression} words'' $\widetilde{z^0}, \widetilde{z^1},
\widetilde{z^2}, \ldots, \widetilde{z^{n^c}}$ such that for $0
\leq s \leq n^c$ there is an equality $\widetilde{z^s}=z^s$ in
$G$. Moreover, we can transform each $z \widetilde {z^s}$ to
$\widetilde{z^{s+1}}$ via a $\PP$-sequence, which when all are
concatenated: $$z^{n^c} \to z^{n^c} \, \widetilde{z^0} \to
z^{n^c-1} \, \widetilde{z^1} \to z^{n^c-2} \, \widetilde{z^2} \to
\cdots \to z \, \widetilde{z^{n^c-1}} \to \widetilde{z^{n^c}},$$
gives a $\PP$-sequence that converts $z^{n^c}$ to
$\widetilde{z^{n^c}}$, and has area at most $\xi n^{c+1}$ and
filling length at most $\xi n$.
\end{cor}

\ni \emph{Proof.} First observe that if we can prove the corollary
for any one finite set of relations $\RR$ then it holds true for
any finite set of relations.

Since the free reduction of $\widetilde{z^0}$ is the empty word,
we can transform  $z^{n^c}$ to $z^{n^c} \, \widetilde{z^0}$ by
free expansion moves.

By definition, each $z \in \AA_c$ is a commutator of length $c$.
Take any $\RR$ such that $\PP = \langle \AA \mid \RR \rangle$ is a
finite presentation for $G$ and  $\RR$ includes all commutators
involving $\cal{A}^{\pm1}$ of length $c+1$. So the presentation
for a free nilpotent group used in Proposition~\ref{compressing
commutators} is a subpresentation of $\PP$ and we can invoke that
proposition to give us the compression words $\widetilde {z^s}$
and the $\PP$-sequence from $z\widetilde {z^s}$ to
$\widetilde{z^{s+1}}$, for each $s$.

The number of $s \in \set{0,1,\ldots, n^c-1}$ such that $n^k \mid
(s+1)$ is at most $n^{c-k}$.   For each such $s$ the area of the
sequence from $z \widetilde{z^s}$ to $\widetilde{z^{s+1}}$ is at
most $\kappa n^{k+1}$ by Proposition~\ref{compressing
commutators}.  So, taking $\xi := (c+1) \kappa$, the total area of
the $\PP$-sequence in question is at most $\xi n^{c+1}$. \qed

\bs

Before we come to the proof of Theorem~\ref{AFL-pair thm} we give
two technical lemmas.

\begin{lemma}~\label{gen set lemma}
Suppose that the abelian group $H$ is generated by the finite set
$\YY$.  Then there is a subset of $\YY$ that freely generates a
free abelian group having finite index in $H$.
\end{lemma}

\ni \emph{Proof.}  This is by induction on $\abs{\YY}$.  There is
nothing to prove if $\abs{\YY}=0$.  If $\abs{\YY} >0$, let $\YY =
\YY_1 \cup \set{y}$ with $\abs{\YY_1} = \abs{\YY}-1$, and let
$\YY_1$ generate $K$.  By inductive hypothesis, $\YY_1$ has a
subset $\YY_2$ that freely generates a free abelian subgroup $L$
with $\abs{ K : L }$ finite.  If $\abs{ H : L }$ is finite then we
are done. Otherwise $H/K$ is infinite cyclic and generated by $y
K$, and then $\YY_2 \cup \set{y}$ has the required property. \qed

\bs

\begin{lemma}~\label{shortening}
Suppose the finite group $H$ of order $t$ is generated by $\YY$.
Then any word $ w \in (\YY \cup \YY^{-1})^{\star}$ with $\ell(w)
\geq t$ contains a subword $u$ with $u = 1$ in $H$.
\end{lemma}

\ni \emph{Proof.}  Let $w_i$ be the prefix of $w$ of length $i$,
for $0 \leq i \leq \ell(w)$.  Then, if $\ell(w) \geq t$, there
exist $i,j$ with $0 \leq i < j \leq \ell(w)$ and $w_i = w_j$ in
$H$, and then $w_j=w_iv$, where $v$ is the required subword of
$w$. \qed

\bs We now give the proof of Theorem~\ref{AFL-pair thm}.  The
following argument is expressed in terms of null-$\PP$-sequences,
however the reader familiar with van~Kampen diagrams will be able
to reinterpret the ideas in geometric terms.

\bs \ni \emph{Proof of Theorem~\ref{AFL-pair thm}.}  It suffices
to find any one finite presentation $\PP$ for $G$ which admits the
required $(\Area, \FL)$-pair.

Let $\AA$ be the generating set for $G$ defined both in \S1 and
before Corollary~\ref{area}. Define $\bar{G}:=
G/{\Gamma_c}$, which is a finitely generated nilpotent group of
class $c-1$.

Suppose $\PP = \langle {\AA} \mid {\RR} \rangle$ is a finite
presentation for $G$. Then define $\bar{\PP}= \langle \bar{\AA}
\mid \bar{\RR} \rangle$ to be the finite presentation for $\bar{G}$ in
which $\bar{\AA}$ is obtained from $\AA$ by removing all generators of
weight $c$, and $\bar{\RR}$ is obtained from $\RR$ by deleting all
occurrences of these generators from the words in $\RR$.

We prove that $(n^{c+1},n)$ is an $(\Area, \FL)$-pair for $G$ by
induction on $c$. The base case $c=1$ concerns abelian groups. The
process of moving letters through words, collecting together and
cancelling every instance of each particular generator, gives the
required $(\Area, \FL)$-pair of $(n^2,n)$.

We now prove the induction step.  Suppose $w$ is a null-homotopic
word in a finite presentation $\PP$ for $G$. Define $n:= \ell(w)$.
Let $\bar{w}$ be the word obtained from $w$ by removing all
occurrences of generators in $\AA_c$. By induction hypothesis,
$\bar{\PP}$ admits $(n^c,n)$ as an $(\Area,\FL)$-pair. Thus there
is a constant $\bar{\lambda}$ depending only on $\bar{\PP}$ such
that there is a null-$\bar{\PP}$-sequence $\bar{S}$:
$$\bar{w}=\bar{w}_0,\bar{w}_1, \ldots, \bar{w}_{\bar{m}} =1$$ for
$\bar{w}$ with $\Area(\bar{S}) \leq \bar{\lambda} n^c$ and
$\FL(\bar{S}) \leq \bar{\lambda} n$. Our intention is to produce a
null-$\PP$-sequence $S$ for $w$ from $\bar{S}$ that demonstrates
that $G$ admits $(n^{c+1},n)$ as an $(\Area,\FL)$-pair.

Now $\AA_c$ generates the abelian group $\Gamma_c$, and it follows
from Lemma~\ref{gen set lemma} above, that we can write $\AA_c =
\AA_{c1} \cup \AA_{c2}$ where $\AA_{c2}$ freely generates a free
abelian group $K$ having finite index $t$ in $\Gamma_c$. Let
$\AA_{c2} =  \set{ z_1, z_2, \ldots, z_k }$.  For $q = 0, 1,
\ldots , n^c  $ let $\widetilde{{z_j}^{q}}$ denote the
\emph{compression word} for ${z_j}^{q}$ of Corollary~\ref{area}.
Extend this sequence of \emph{compression words} to all
non-negative powers of $z_j$ by defining
$$\widetilde{{z_j}^{A+Bn^c}} \ := \ \widetilde{{z_j}^{A}} \,
(\widetilde{{z_j}^{n^c}})^B$$ for $0 \leq A \leq n^c -1$ and $B
\geq 0$.

The elements of $\AA_c$ are all central in $G$, so we may as well
assume that $\RR$ includes the commutators $[x,z]$ for all $x \in
\AA ^{\pm 1}$ and $z \in {\AA_c}^{\pm 1}$.

Lemma~\ref{app of a relator} tells us that, by suitably adjusting
the constant $\bar{\lambda}$, we can assume that each time a word
$\bar{w}_{i+1}$ is obtained from a word $\bar{w}_i$ in the
sequence $\bar{S}$ by \emph{applying a relator} $\bar r  \in \bar
\RR$, the whole of some cyclic conjugate of $\bar r^{\pm 1}$ is
inserted into $\bar{w}_i$.

Suppose a word $u$ of length at most $t$ in the letters
${\AA_{c1}}^{\pm 1}$ represents an element of $K$. Then we can
choose a word $v_u$ in the letters ${\AA_{c2}}^{\pm 1}$ such that
$u=v_u$ in $G$.  There are only finitely many words of length at
most $t$ so we may as well assume that each $u{v_u}\inv$ is in
$\RR$.

We now explain how to find a null-$\PP$-sequence for $w$, by
\emph{expanding} $\bar S$.  Suppose that the move from
$\bar{w}_{i}$ to $\bar{w}_{i+1}$ is an instance of an
\emph{application of a relator} move in the
null-$\bar{\PP}$-sequence $\bar{S}$.  So $\bar{w}_i = \alpha
\beta$ for some words $\alpha$ and $\beta$.  Also $\bar{w}_{i+1} =
\alpha \bar{r} \beta$, where some cyclic conjugate of $\bar{r}$ is
in $\bar{\RR}^{\pm 1}$.  There is some $r$ which has a cyclic
conjugate in $\RR^{\pm 1}$ and from which one can obtain $\bar r$
by removing all occurrences of generators in ${\AA_c}^{\pm 1}$.
The idea is to apply an \emph{expanded version} of the move:
instead of inserting $\bar{r}$, we insert $r$.  The newly
introduced generators from ${\AA_c}^{\pm 1}$ are then collected at
the ends of the word and \emph{compressed} before the
$\PP$-sequence is continued as dictated by $\bar{\PP}$.

More precisely, we construct a null-$\PP$-sequence $S$ for $w$ as
follows.  Amongst the words in $S$ will be words of the form $$w_i
\ = \ (\widetilde{{z_k}^{q_{ki}^-}} )^{-1} \ldots \,
(\widetilde{{z_2}^{q_{2i}^-}})^{-1} \,
(\widetilde{{z_1}^{q_{1i}^-}})^{-1} \, \bar{w}_{i} \, u_i \,
\widetilde{{z_1}^{q_{1i}^+}} \, \widetilde{{z_2}^{q_{2i}^+}}\ldots
\widetilde{{z_k}^{q_{ki}^+}},$$ for $0 \leq i \leq \bar{m}$.  In
these words each $q_{ji}^+$ and $q_{ji}^-$ is a non-negative
integer, and $u_i$ is a word of length less than $t$ in the
generators ${\AA_{c1}}^{\pm 1}$.

The $\PP$-sequence $S$ starts by transforming $w$ to $w_0$ by
moving all generators from ${\AA_c}^{\pm 1}$, through $w_i$ using
commutator relations (recall that we assumed these to be in
$\RR$): move all occurrences of $z_1, z_2, \ldots, z_k$ towards
the (right-hand) end of the word, all occurrences of ${z_1}^{-1},
{z_2}^{-1}, \ldots, {z_k}^{-1}$ towards the start of the word, and
all generators from ${\AA_{c1}}^{\pm 1}$ towards the right to form
the subword $u_0$.  These generators from $\AA_c$ are all
\emph{collected} and \emph{compressed} as we now explain.

When a letter $z_j$ is moved towards the (right-hand) end of a
word it arrives at some subword word $\widetilde{{z_j}^{q}}$,
where $q \geq 0$. We then transform $$z_j \, \widetilde{{z_j}^{q}}
\ \to \ \widetilde{{z_j}^{q+1}}$$ via a $\PP$-sequence $S_{jq}$ of
Corollary~\ref{area}. Similarly when a ${z_j}^{-1}$ is moved
towards the start of the word it meets some
$\left(\widetilde{{z_j}^{q}}\right)^{-1}$, where $q \geq 0$, and
we then make the transformation
$$\left(\widetilde{{z_j}^{q}}\right)^{-1} \, {z_j}^{-1} \ \to \
\left(\widetilde{{z_j}^{q+1}}\right)^{-1}$$ using the
$\PP$-sequence obtained by inverting every word in $S_{jq}$.

The letters from ${\AA_{c1}}^{\pm 1}$ are collected together at
the appropriate position in $w_0$.  However in the process of
collection, if we create a subword that represents an element of
$K$ we immediately replace this subword by a word in
${\AA_{c2}}^{\pm 1}$ (which, by assumption, we can do using a
relator in $\RR$). The resulting letters from ${\AA_{c2}}^{\pm 1}$
are then moved and collected in their compression words in the
manner already explained above. It follows from
Lemma~\ref{shortening} that the collected words in the letters
${\AA_{c1}}^{\pm 1}$ always has length less than $t$.

The $\PP$-sequence $S$ continues with a concatenation of
$\PP$-sequences that transform $w_i$ to $w_{i+1}$ for $i = 0,1,
\ldots, \bar{m}-1$. The sequence from $w_i$ to $w_{i+1}$ starts by
transforming the subword $\bar{w}_i$ as dictated by $\bar{S}$: if
the move is a free expansion or a free reduction then this move is
performed and we immediately have $w_{i+1}$. However if the move
is an \emph{application of a relator} then we apply the
\emph{expanded version}.  This introduces letters from
${\AA_c}^{\pm 1}$. Move the $z_1, z_2, \ldots, z_k$ towards the
end of the word, the occurrences of ${z_1}^{-1}, {z_2}^{-1},
\ldots, {z_k}^{-1}$ towards the start of the word, and all
generators from ${\AA_{c1}}^{\pm 1}$ to the start of the subword
$u_i$, and \emph{compress} in the manner explained above.

\ms Now $$w_{\bar{m}} \ = \ (\widetilde{{z_k}^{q_{k{\bar{m}}}^-}}
)^{-1} \ldots \, (\widetilde{{z_2}^{q_{2{\bar{m}}}^-}})^{-1} \,
(\widetilde{{z_1}^{q_{1{\bar{m}}}^-}})^{-1} \,  u_{\bar{m}} \,
\widetilde{{z_1}^{q_{1{\bar{m}}}^+}} \,
\widetilde{{z_2}^{q_{2{\bar{m}}}^+}}\ldots
\widetilde{{z_k}^{q_{k{\bar{m}}}^+}},$$ is a null-homotopic word
in $\Gamma_c$, and so in fact $u_{\bar{m}}$ must be the empty word
(else it would represent an element of $K$, and hence it would
have already been eliminated from $w_{\bar m}$ by the procedure
above).  So $q_{j{\bar{m}}}^+ = q_{j{\bar{m}}}^-$ for $j=1,2,
\ldots, k$ and we can complete our $\PP$-sequence for $w$ by
freely reducing $w_{\bar m}$ to the empty word.

\ms All that remains is to explain why there is a constant
$\lambda$, depending only on $\PP$, such that
\begin{eqnarray*}
\Area(S) & \leq & \lambda \, n^{c+1} \\ %
\FL(S) & \leq & \lambda \, n.
\end{eqnarray*}

Let $M$ be the maximum number of letters from ${\AA_c}^{\pm 1}$
occurring in any one of the relators in $\RR$.  At most $M
\Area(\bar{S})$ letters from ${\AA_{c1}}^{\pm 1}$ are collected
and compressed in the subwords $u_i$.  In the compression process
one relator may be used for each of these letters, in the process
releasing further letters from ${\AA_{c2}}^{\pm 1}$ to be
collected.

Letters from ${\AA_{c2}}^{\pm 1}$ are released in $S$ on account
of each application of a relator in $\bar S$ and also each relator
used when compressing the word that collects letters from
${\AA_{c1}}^{\pm 1}$. Thus the powers $q_{ji}^+$ and $q_{ji}^-$
are all bounded by $2M \Area(\bar{S})$.

Corollary~\ref{area} therefore tells us that the length of the
compression words for powers of elements of $\AA_{c2}$ is $O(n)$.
Moreover the total contribution to the area of $S$ arising in the
compression process is $O(n^{c+1})$.

It follows that there is a constant $\lambda'>0$, depending only
on $\PP$, such that $\FL(S) \leq \lambda' \, n$.  For the bound on
$\Area(S)$ notice that the observations above tell us that the
total contributions on account of the compression process is
$O(n^{c+1})$, and we need only bound the contribution made by our
use of commutators to move letters through words.  But each of the
$O(n^c)$ letters from ${\AA_c}^{\pm 1}$ that arise in $S$ is moved
a distance at most $\lambda' \, n$.  So indeed, we can find
$\lambda \geq \lambda'$, depending only on $\PP$, such that
$\Area(S) \leq \lambda \, n^{c+1}$ and $\FL(S) \leq \lambda \, n$.
\qed

\bibliographystyle{plain}
\bibliography{bibli}

\small{ \ni \textsc{Steve M.\ Gersten} \rule{0mm}{6mm} \\
Mathematics Department, 155S. 1400E., Rm. 233, Univ. of Utah, Salt
Lake City, UT 84112, USA
\\ \texttt{gersten@math.utah.edu, \
http:/\!/www.math.utah.edu/$\sim$gersten}

\ni \textsc{Derek F.\ Holt} \rule{0mm}{6mm} \\ Mathematics
Institute, University of Warwick, Coventry, CV4 7AL, UK \\
\texttt{dfh@maths.warwick.ac.uk, \
http:/\!/www.maths.warwick.ac.uk/$\sim$dfh}

\ni  \textsc{Tim R.\ Riley} \rule{0mm}{6mm} \\ Mathematical
Institute, 24-28 St. Giles', Oxford OX1 3LB, UK \\
\texttt{rileyt@maths.ox.ac.uk, \
http:/\!/www.maths.ox.ac.uk/$\sim$rileyt } }
\end{document}